\documentclass[11pt,twosided]{article}
\usepackage{dsfont}
\usepackage{bm}
\usepackage{float}
\usepackage{srcltx}
\usepackage{mathrsfs}
\usepackage{graphicx}
\usepackage{amsfonts}
\usepackage{indentfirst}
\usepackage{amssymb,amsfonts,amsmath, psfrag,
}
\usepackage{diagbox}
\usepackage{caption}

\usepackage{leftidx}
\usepackage{pgf,tikz}
\usepackage{cite}
\usetikzlibrary{shapes,arrows}
\allowdisplaybreaks[4]

\usepackage{color}

\usepackage{amssymb,color}

\definecolor{c20}{rgb}{0.,0.7,0.}
\definecolor{c30}{rgb}{0.,0.,1.}
\definecolor{c40}{rgb}{1,0.1,0.7}
\definecolor{c50}{rgb}{1,0,0}
\definecolor{c60}{rgb}{0,0.9,0.1}

\usepackage{amssymb,color}

\newcommand{\BQN}{\begin{eqnarray}}
\newcommand{\EQN}{\end{eqnarray}}
\newcommand{\BQNY}{\begin{eqnarray*}}
\newcommand{\EQNY}{\end{eqnarray*}}

\newcommand{\BS}{\begin{sat}}
\newcommand{\ES}{\end{sat}}
\newcommand{\BT}{\begin{theo}}
\newcommand{\ET}{\end{theo}}
\newcommand{\BK}{\begin{korr}}
\newcommand{\EK}{\end{korr}}

\newcommand{\BD}{\begin{de}}
\newcommand{\ED}{\end{de}}

\newcommand{\BIT}{\begin{itemize}}
\newcommand{\EIT}{\end{itemize}}
\newcommand{\BDI}{\begin{description}}
\newcommand{\EDI}{\end{description}}

\newcommand{\BRM}{\begin{remarks}}
\newcommand{\ERM}{\end{remarks}}

\newcommand{\BEL}{\begin{lem}}
\newcommand{\EEL}{\end{lem}}

\newtheorem{theo}{Theorem}[section]

\newtheorem{lem}[theo]{Lemma}

\newtheorem{coro}[theo]{Corollary}

\newtheorem{remarks}[theo]{Remarks}

\newcommand{\COM}[1]{}

\begin{document}
\title{The joint distribution of the Parisian ruin time and the number of claims until Parisian ruin in the classical risk model}
\author{Irmina Czarna\thanks{Mathematical Institute, University of Wroclaw, Poland. E-mail: czarna@math.uni.wroc.pl}\ \ Yanhong Li\thanks{School of mathematical Sciences and LPMC Nankai University, Tianjin 300071, P.R. China.
E-mail: yanhongli@mail.nankai.edu.cn}\ \ \ Zbigniew Palmowski\thanks{Mathematical Institute, University of Wroclaw, Poland. E-mail: zbigniew.palmowski@gmail.com} \ \  Chunming Zhao\thanks{Department of Statistics
School of Mathematics Southwest Jiaotong University, Chengdu, Sichuan, 611756, PR of China. Email:cmzhao@swjtu.cn}}
\date{}
\maketitle
\begin{abstract}
In this paper we propose new iterative algorithm of calculating the joint distribution of the Parisian ruin time and the number of claims until Parisian ruin for the classical risk model.
Examples are provided when the generic claim size is exponentially distributed.
\end{abstract}


\section{Introduction}
The distribution of the number of claims until ruin has been
in the centre of interest for many years. One of the first references dealing with this problem
is Beard \cite{Beard}. The main first step was done by Stanford and Stroi\'{n}ski \cite{Stanford}
who produced recursive procedures to calculate the probability of ruin at
the $n$th claim arrival epoch in the classical risk model.
Eg\'{i}do dos Reis \cite{Egido} derived the moment generating function of the number of claims
until ruin in the classical risk model. He inverted this for certain claim size
distributions, and, using a duality argument, found moments of the number
of claims until ruin when the initial surplus is $0$.
The next main step was done
by Landriault et al. \cite{LaShWi2011} who considered a Sparre Andersen risk
model with exponential claims. Using Gerber-Shiu type analysis (see Gerber
and Shiu \cite{GeSh2005})) they derived a number of results including an expression
for the probability function of the number of claims until ruin. The main idea of getting these nice results
followed approach of Dickson and Willmot \cite{DiWi2005}.
The main results of our paper are closely related with
the seminal paper of Dickson \cite{Di2012} who using probabilistic arguments
derived the expression for the joint
density of the time of ruin and the number of claims until ruin in the
classical risk model. From this he obtained a general expression for
the probability function of the number of claims until ruin. He also
considered the moments of the number of claims until ruin and illustrate
all results in the case of exponentially distributed individual claims.
Frostig et al.\cite{FrPiPo2012} and Zhao and Zhang \cite{ZhZh2013} analyzed similar problems.

In this paper we extend results concerning classical ruin into so-called Parisian type of ruin.
This type of ruin occurs if the surplus process falls below zero and stays below zero for a continuous time of interval of length $d$; see Figure 1.
We believe that the Parisian ruin probability and other related quantities
might be more appropriate measures of risk than the ones identified for the classical ruin.
The main reason is that it gives the insurance companies the chance to achieve solvency.
The idea of Parisian ruin comes from Parisian options which  was first introduced by Chesney et al. \cite{ChJeYo1997}. Dassios and Wu \cite{DaWu2008} considered the Parisian ruin probability for the classical risk model with exponential claims and for the Brownian motion with drift. Czarna and Palmowski \cite{CzPa2011} and Loeffen et al. \cite{LoCzPa2013} analyzed the Parisian ruin probability for a general spectrally negative L\'evy process. Other relevant papers are Landriault et al. \cite{LaReZh2011} and \cite{LaReZh2014}, where the deterministic and fix delay $d$ is replaced by an independent exponential random variable.

Our paper in a sense has similar goal like in Dickson \cite{Di2012} and Landriault et al. \cite{LaReZh2011}, that is we want
to identify the joint density of the time of Parisian ruin and the number of claims until Parisian ruin.
Although our focus is more consistent with Stanford and Stroi\'{n}ski \cite{Stanford} - we want
to create efficient iterative algorithm of finding above quantity.

Formally, in this paper we consider a continuous-time surplus process:
\begin{equation}\label{model}
   U(t):=u+ct-\sum_{i=1}^{N_t}X_{i},
\end{equation}
where the non-negative constant $u $ denotes the initial reserve, the positive constant $c$ is the rate of premium income, $N_t$ describes the number of claims counted up to time $t$ which is a Poisson process with parameter $\lambda $ and $\{X_{i}\}_{i=1}^{\infty}$ are claim sizes which are independent and identically distributed non-negative  random variables that are also independent of $N_t$. We denote by $F(x)$ and  $f(x)$ the distribution function and density function, respectively. We assume  $c>\lambda\mathbb{E}(X_{1})$ assuring that  ruin is not certain.
\begin{center}
\definecolor{ffqqtt}{rgb}{1,0,0.2}
\begin{tikzpicture}[line cap=round,line join=round,>=triangle 45,x=1.0cm,y=1.0cm]
\draw[->,color=black] (-1,0) -- (15,0);
\foreach \x in {}
\draw[shift={(\x,0)},color=black] (0pt,-2pt);
\draw[->,color=black] (0,-3.5) -- (0,4);
\foreach \y in {}
\draw[shift={(0,\y)},color=black] (2pt,0pt) -- (-2pt,0pt);
\clip(-1,-3,5) rectangle (15,4);
\draw (0,1)-- (1,2);
\draw [dash pattern=on 3pt off 3pt] (1,2)-- (1,0.5);
\draw (1,0.5)-- (2,1.5);
\draw [dash pattern=on 3pt off 3pt] (2,1.5)-- (2,-1);
\draw (2,-1)-- (5,2);
\draw [dash pattern=on 3pt off 3pt] (5,2)-- (5,1);
\draw (5,1)-- (7,3);
\draw [dash pattern=on 3pt off 3pt] (7,3)-- (7,-3);
\draw (7,-3)-- (9,-1);
\draw [dash pattern=on 3pt off 3pt] (9,-1)-- (9,-2);
\draw (9,-2)-- (10.5,-0.5);
\draw [dash pattern=on 3pt off 3pt] (10.5,-0.5)-- (10.5,-1);
\draw (10.5,-1)-- (12,0.5);
\draw [shift={(3,-3.26)},color=ffqqtt]  plot[domain=1.27:1.87,variable=\t]({1*3.41*cos(\t r)+0*3.41*sin(\t r)},{0*3.41*cos(\t r)+1*3.41*sin(\t r)});
\draw [shift={(8,-3.26)},color=ffqqtt]  plot[domain=1.27:1.87,variable=\t]({1*3.41*cos(\t r)+0*3.41*sin(\t r)},{0*3.41*cos(\t r)+1*3.41*sin(\t r)});
\draw (-0.95,4) node[anchor=north west] {$U(t)$};
\draw (14.5,0.07) node[anchor=north west] {$t$};
\draw (-0.56,0.14) node[anchor=north west] {$0$};
\draw (-0.53,1.36) node[anchor=north west] {$u$};
\draw (2.82,0.79) node[anchor=north west] {$d$};
\draw (8.04,0.84) node[anchor=north west] {$d$};
\draw (8.86,0.15) node[anchor=north west] {$\tau_u^d$};
\draw (1.9,0) node[anchor=north west] {$\tau_u$};
\draw (0.5,-3.51) node[anchor=north west] {${\rm Fig. 1. \; Surplus  \; process}  \; U(t) \; {\rm and \; Parisian \; ruin.}$};
\end{tikzpicture}
\end{center}

 We define the Parisian time of ruin by
 \begin{equation*}
 \tau_u^{d}:=\inf\{t>0: t-\sup\{s<t: U(s)\ge0\}\ge d, U(t)<0\}.
 \end{equation*}
 We denote the joint density of $N_{\tau_u^d}$ and $\tau_u^d$ by (hereafter $\mathbb{N}=\{0,1,2,\cdots\}$):
 \begin{equation*}
 w_u^{d}(n,t):=\frac{d}{dt}\psi_u^d(n,t),\ \ \ n\in\mathbb{N}, t\geq 0
 \end{equation*}
 with
\begin{equation}\label{pdnt}
\psi_u^d(n,t):=\mathbb{P}(N_{\tau_u^d}=n, \tau_u^d\leq t|U(0)=u), \ \ \ n\in\mathbb{N}, t\geq d.
 \end{equation}
 Further, let $p_u^d(n)$ denotes the probability that there have been exactly $n$ claims up to Parisian ruin event, so that
 \begin{eqnarray}
 p_u^d(n)&:=&\mathbb{P}(N_{\tau_u^d}=n, \tau_u^d <\infty|U(0)=u)=\int_d^{\infty}w_u^{d}(n,t)dt. \label{pdn}
 \end{eqnarray}
The main goal of this paper is to give efficient iterative algorithm of calculating of
$w_u^{d}(n,t)$ and hence $p_u^d(n)$.

Note that the $d=0$ case corresponds to the classical ruin problem.
Then we deal with the classical ruin time of the risk process (\ref{model}):
\begin{equation}\label{classruin}
 \tau_u:=\inf\{t\ge0: U(t)<0\}.
\end{equation}

The rest of the paper is organized as follows. In Section 2
we give the main representations of the joint density of $N_{\tau^d}$ and $\tau^d$ and
prove the main results. Finally, in Section 3 we analyze some particular examples
and give extensive numerical analysis.
\textbf{}

\section{Main representation}

Since Parisian ruin occurs if the surplus falls below zero and stays below zero for a continuous time interval of length $d$,
Parisian ruin time must be larger than $d$. Throughout this paper, we can assume that $w_u^d(n,t)=0$ when $t\leq d, u\ge0$.
For $t>d$ and $u\geq 0$ we will now identify the joint density $w_u^d(n,t)$ ($n\in\mathbb{N}$)
of $N_{\tau_u^d}$ and $\tau_u^d$.

\subsection{The expression of the joint density $w_u(k,t,y)$}\label{secDickson}

Our results heavily use the main result of Dickson \cite{Di2012}.
He considers the joint density of the number $N_{\tau_u}$
of claims until ruin (including the ruin-caused claim), $\tau_u$ given in (\ref{classruin})
and the deficit at ruin $|U(\tau_u)|$
defined by :
$$
w_u(k,t,y):=\frac{\partial^2}{\partial t\partial y}\psi_u(k,t,y),\ \ \ k\in\mathbb{N}, t\ge0, y>0
$$
with
$$\psi_u(k,t,y):=\mathbb{P}(N_{\tau_u}=k, \tau_u\leq t, |U(\tau_u)|\leq y \big| U(0)=u), \ \ \ k\in\mathbb{N}, t\ge0, y>0.$$

For any function $g$ we denote by $g^{k*}$ ($k\geq 0$) the $k$-fold convolution of $g$ with itself,
where $g^{1*}=g$ and $g^{0*}(t)=\delta_0(t)$ for the impulse function $\delta_0$ at $0$.
We are now in position to state the main result of Dickson \cite{Di2012}.
\begin{theo}
For $k=1,2,\ldots$  we have that
\begin{eqnarray}\label{exam_1_inverse_laplace}
   w_0(k,t,y)&=&\int_{0}^{ct}\frac{x}{ct}e^{-\lambda t}\frac{\lambda^k t^{k-1}}{(k-1)!}f^{(k-1)*}(ct-x)f(x+y)dx
   \end{eqnarray}
and for $u>0$
\begin{eqnarray*}
   w_u(k,t,y)&=&\int_0^{u+ct}e^{-\lambda t}\frac{\lambda^kt^{k-1}}{(k-1)!}f^{(k-1)*}(u+ct-x)f(x+y)dx\\
   &&-c\sum_{j=1}^{k-1}\int_0^t e^{-\lambda s}\frac{(\lambda s)^j}{j!}f^{j*}(u+cs)\omega_0(k-j,t-s,y)ds.
\end{eqnarray*}
\end{theo}
\begin{coro}
When the generic claim amount is exponentially distributed with mean $1/\mu$, i.e. $f(x)=\mu e^{-\mu x}$, then
\begin{eqnarray}
&&w_0{(k,t,y)} =\frac{\lambda^k \mu^{k}c^{k-1}}{k!(k-1)!}t^{2k-2}e^{-(\lambda+\mu c)t}  e^{-\mu y},\label{w0kt}\\
&&w_{u}{(k,t,y)}=\frac{\lambda^k \mu^{k}(ku+ct)(u+ct)^{k-2}}{k!(k-1)!}t^{k-1}e^{-(\lambda+\mu c)t-\mu u} e^{-\mu y}.\\
\end{eqnarray}
\end{coro}
For $u\geq 0$ we denote by the density of having $k$ claims up to the classical ruin time $\tau_u$ and having the deficit $y$ at the ruin:
\begin{eqnarray*}
 w_{u}(k,y)&=&\int_0^{\infty}w_{u}(k,t,y)dt.
\end{eqnarray*}
In particular, for exponential claim size with intensity $\mu$ we have:
\[w_0(k,y)=\frac{(2k-2)!}{k!(k-1)!}\frac{\lambda^k \mu^{k}c^{k-1}}{(\lambda+\mu c)^{2k-1}}  e^{-\mu y}.\]

 \subsection{The joint distribution of the first upward passage time and the number of claims}
 For $y\ge0$ we define the first upward passage time of our classical risk process:
 \begin{equation*}
 \tau_y^+=\min\{s\ge0,U(s)=y|U(0)=0\}.
 \end{equation*}
 We denote by
 \begin{equation*}
 v_y(k,t):=\frac{d}{dt}V_y(k,t),\ \ \ k\in\mathbb{N}, t\ge y/c
 \end{equation*}
 the density of having $k$ jumps up to first passage time of level $y$ that happens at time $t$, that is:
\begin{equation*}
V_y(k,t):=\mathbb{P}(N_{\tau_y^+}=k, \tau_y^+\leq t|U(0)=0), \ \ \ k\in\mathbb{N}, t\ge 0.
 \end{equation*}

\begin{theo}
We have:
\begin{eqnarray}\label{v y k t}
v_y(k,t)&=&\frac{\lambda^k}{k!}yt^{k-1}e^{-\lambda t}f^{k*}(ct-y).
\end{eqnarray}
\end{theo}
\textbf{Proof.}
For $r\in(0,1]$ and $\delta>0$, we define the bivariate Laplace transform of $(\tau_y^+, N_{\tau_y^+})$:
\begin{eqnarray}\label{define phi}
    \phi(y)&:=&E[r^{N_{\tau_y^+}}e^{-\delta\tau_y^+}\mathbb{I}(\tau_y^+<\infty)|U(0)=0]\nonumber\\
           &=&\int_0^\infty e^{-\delta t}\sum_{k=0}^\infty r^k v_y(k,t)dt.
\end{eqnarray}

Considering an infinitesimal time interval $(0,dt)$ we have:
\begin{eqnarray}\label{phi}
\phi(y)&=&e^{-\delta dt}e^{-\lambda dt}\phi(y-cdt)+re^{-\delta dt}(1-e^{-\lambda dt})\int_0^{\infty}\phi(y-cdt+x)f(x)dx+o(dt)\nonumber\\
&=&[1-(\lambda+\delta)dt]\phi(y-cdt)+\lambda rdt\int_0^{\infty}\phi(y-cdt+x)f(x)dx+o(dt).
\end{eqnarray}
Subtracting $\phi(y-ct)$ from both sides of above equation, multiplying by $1/dt$ and letting $dt\rightarrow 0$ produce the following integro-differential equation:
\begin{equation}\label{integro differential equation}
c\phi^{'}(y)=-(\lambda+\delta)\phi(y)+\lambda r\int_0^{\infty}\phi(y+x)f(x)dx.\end{equation}
Clearly, when $y=0$, we have
\begin{equation}\label{boundary condition}
\phi(0)=1.
\end{equation}
Since the solution to (\ref{integro differential equation}) with boundary condition (\ref{boundary condition}) is unique, we assume that $\phi(y)$ is of the form
\begin{equation*}
\phi(y)=c(y)e^{-by}.
\end{equation*}
The boundary condition $\phi(0)=1$ gives $c(y)=1$, so that $\phi(y)=e^{-by}$. Note that the real part of $b$ must be positive, because otherwise it would be a contradiction to the fact that $\lim_{y\rightarrow\infty}\phi(y)=0$. It is known that the Lundberg's fundamental equation of the classical risk model is given by
\begin{equation*}
\lambda+\delta-cs=\lambda r\hat f (s).
\end{equation*}
We denote the positive solution by $\rho$.

Now, using a similar approach as in Li \cite{Li2008a}, Zhao and Zhang \cite{ZhZh2013} we can
obtain the solution of the integro-differential equation (\ref{integro differential equation}):
\begin{equation*}
\phi(y)=e^{-\rho y}.
\end{equation*}

We recall now the \textbf{Lagrange's Expansion Theorem} (see Lagrange \cite[p. 251-326]{La1770}).
Given two functions $\alpha(z)$ and $\beta(z)$ which are both analytic on and inside a contour $D$ surrounding a point $a$, if $r$ satisfies the inequality
\begin{equation}
|r \beta(z)|<|z-a|,
\end{equation}
for every $z$ on the perimeter of $D$, then $z-a-r \varphi(z)$, as a function of $z$, has exactly one zero $\eta$ in the interior of $D$,  and
we have further
\begin{equation}
    \alpha(\eta)=\alpha(a)+\sum_{k=1}^{\infty}\frac{r^k}{k!}\frac{d^{k-1}}{dx^{k-1}}\big(\alpha'(x)\beta^k(x)\big)\!\big|_{x=a}.
\end{equation}

It follows from above fact that:
\begin{equation}\label{e rho y}
e^{-\rho y}=\int_0^{\infty}e^{-\delta t}\sum_{k=0}^{\infty}r^k\left(\frac{\lambda^k}{k!}yt^{k-1}e^{-\lambda t}f^{k*}(ct-y)\right)dt.
\end{equation}
Comparing (\ref{define phi}) and (\ref{e rho y}) gives the assertion of the theorem.
\vspace{3mm} \hfill $\Box$

\begin{coro}
When the individual claim amounts are exponentially distributed with mean $1/\mu$ then
\begin{eqnarray}\label{expv}
v_y(k,t)=\left\{
  \begin{array}{lll}
&\frac{y}{t}e^{-\lambda t}\delta_0(ct-y),&k=0,\\
&\frac{\lambda^k \mu^{k}y}{k!(k-1)!}t^{k-1}(ct-y)^{k-1}e^{-(\lambda+\mu c)t}  e^{\mu y},&k>0.
  \end{array}
\right.
\end{eqnarray}
\end{coro}

\subsection{The expression of $w_u^d(n,t)$}
Recall that $w_u^d(n,t)$ is the joint density that Parisian ruin occurs at time $t$ and  there are $n$ claims up to time $t$.
The main result of this paper gives recursive algorithm of calculating
the density $w_u^d(k,t)$.
\begin{theo}\label{maintheorem}
We have
\begin{equation}\label{w u d 1}
w_u^d(1,t)=\lambda e^{-\lambda t}\bar{F}(u+ct)
\end{equation}
and for $n=2,3,\ldots$,
  \begin{eqnarray}\label{u>0}
  \lefteqn{w_u^d(n,t)=}\\
  &=&\sum_{k=0}^{n-1}\int_{cd}^{\infty}w_{u}(n-k,t-d,y)\frac{(\lambda d)^k}{k!}e^{-\lambda d}dy+\sum_{k=1}^{n-1}\int_0^{cd}w_{u}(n-k,t-d,y)\nonumber\\
   &&\times\Big[\frac{(\lambda d)^k}{k!}e^{-\lambda d}\bar{F}^{k*}(cd-y)-\sum_{m=0}^{k-1}\int_{\frac{y}{c}}^dv_y(m,s) \frac{(\lambda (d-s))^{k-m}}{(k-m)!}e^{-\lambda (d-s)}\bar{F}^{(k-m)*}(c(d-s))ds \Big]dy\nonumber\\
    &&+\sum_{l=1}^{n-1}\sum_{k=0}^{n-l-1}\int_0^{cd}\int_{\frac{y}{c}}^d\int_0^{\max(t-d-s,0)}w_{u}(l,t_1,y)v_y(k,s)w_0^d(n-l-k,t-t_1-s)dt_1dsdy.
     \end{eqnarray}
\end{theo}
\textbf{Proof.}\
Fact that the Parisian ruin occurs in the time interval $(0,t]$ and there is only one claim up to Parisian ruin time, means that
only one claim  occurs before time $t-d$ that cause the classical ruin, the deficit is larger than $cd$ and there will be no claims
within time that risk process spent below $0$, see Figure 2. This gives (\ref{w u d 1}).
\begin{center}
\definecolor{ffttqq}{rgb}{1,0.2,0}
\begin{tikzpicture}[line cap=round,line join=round,>=triangle 45,x=1.0cm,y=1.0cm]
\draw[->,color=black] (-1.04,0) -- (12,0);
\foreach \x in {}
\draw[shift={(\x,0)},color=black] (0pt,-2pt);
\draw[->,color=black] (0,-3.46) -- (0,2.5);
\foreach \y in {}
\draw[shift={(0,\y)},color=black] (2pt,0pt) -- (-2pt,0pt);
\clip(-1.04,-4) rectangle (12,2.5);
\draw (0,1)-- (1,2);
\draw (1,-3)-- (5,1);
\draw [dash pattern=on 3pt off 3pt] (1,2)-- (1,-3);
\draw [shift={(2,-3.26)},color=ffttqq]  plot[domain=1.27:1.87,variable=\t]({1*3.41*cos(\t r)+0*3.41*sin(\t r)},{0*3.41*cos(\t r)+1*3.41*sin(\t r)});
\draw (-0.9,2.5) node[anchor=north west] {$U(t)$};
\draw (-0.42,0) node[anchor=north west] {$0$};
\draw (0.98,0) node[anchor=north west] {$\tau_u$};
\draw (1.86,0.7) node[anchor=north west] {$d$};
\draw (2.86,0.15) node[anchor=north west] {$\tau_u^d$};
\draw (11.5,0.1) node[anchor=north west] {$t$};
\draw (-0.42,1.1) node[anchor=north west] {$u$};
\draw (0.5,-3.3) node[anchor=north west] {{\rm Fig. 2.  The case of Parisian ruin with one claim}};
\end{tikzpicture}
\end{center}

The arguments behind the formula (\ref{u>0}) are as follows.
We know that the Parisian ruin occurs after classical ruin.  There are only two cases:
 \begin{itemize}
\item $\tau^d_u=\tau_u+d$, i.e., the surplus will stay below zero for a continuous time interval of length $d$ after the classical ruin time.  Let us assume that there are $k$ $(0\leq k\leq n-1)$ claims during the interval $(\tau_u,\tau_u+d]$ and $n-k$ claims during the interval $(0,\tau_u]$. If the deficit at the classical ruin is more than $cd$ then the surplus can not exceed $0$ before $\tau_u+d$ no matter how much the cumulative amount of the $k$ claims is. This covers the first term of formula (\ref{u>0}). However, if the deficit is less than $cd$ (formulated as the second term) then it also includes the possibility that the surplus has been up-crossing $0$ prior to time $\tau_u+d$ which should be subtracted. To take into account w suppose that $\tau_u+s$ is the first time before $\tau_u+d$ at which there was an up-crossing of the surplus process through $0$ and there are $m$ claims during the interval $(\tau_u,\tau_u+s]$ and hence $k-m$ claims during the interval $(\tau_u+s,\tau_u+d]$.
 \item $\tau^d_u>\tau_u+d$, i.e., the surplus exceeds $0$ in the interval $(\tau_u,\tau_u+d]$ (we assume also that classical ruin happens at time $t_1$). We apply probabilistic arguments to construct the last term of (\ref{u>0}); see Prabhu\cite{Pr1961}. We take $\tau_u+s$ to be the first time before $\tau_u+d$ at which there
 is an up-crossing of the surplus process through $0$. Further, we suppose that there are $l$ claims in $(0,\tau_u]$ and $k$ claims in $(\tau_u,\tau_u+s]$. Additionally,
 when risk process up-crosses zero it does in continuous way. So we can restart the our considerations with $u=0$, the Parisian ruin time equal to $t-t_1-s$ and
 $(n-l-k)$ amount of claims counted up to this time.
 \end{itemize}
\vspace{3mm} \hfill $\Box$

In particular, from (\ref{u>0}) for $u=0$ we can obtain the following corollary.
\begin{coro}
\begin{eqnarray}\label{u=0}
    \lefteqn{w_0^d(n,t)}\\&=&\sum_{k=0}^{n-1}\int_{cd}^{\infty}w_{0}(n-k,t-d,y)\frac{(\lambda d)^k}{k!}e^{-\lambda d}dy+\sum_{k=1}^{n-1}\int_0^{cd}w_{0}(n-k,t-d,y)\nonumber\\
   &&\times\Big[\frac{(\lambda d)^k}{k!}e^{-\lambda d}\bar{F}^{k*}(cd-y)-\sum_{m=0}^{k-1}\int_{\frac{y}{c}}^dv_y(m,s) \frac{(\lambda (d-s))^{k-m}}{(k-m)!}e^{-\lambda (d-s)} \bar{F}^{(k-m)*}(c(d-s))ds\Big]dy\nonumber\\
    &&+\sum_{l=1}^{n-1}\sum_{k=0}^{n-l-1}\int_0^{cd}\int_{\frac{y}{c}}^d\int_0^{\max(t-d-s,0)}w_{0}(l,t_1,y)v_y(k,s)w_0^d(n-l-k,t-t_1-s)dt_1dsdy.
 \end{eqnarray}
 \end{coro}

 In order to find the explicit expression of $w_0^d(n,t)$, we first consider $w_0^d(1,t)$ which denotes the joint density function when Parisian ruin occurs at time $t$ and there is only one claim up to time $t$.  Plugging $u=0$ into  (\ref{w u d 1}) we have:
 \begin{equation}\label{w 0 d 1}
w_0^d(1,t)=\lambda e^{-\lambda t}\bar{F}(ct),    \ \ \ \ \ t>d.
\end{equation}
Then substituting (\ref{w 0 d 1}), (\ref{v y k t}) and (\ref{exam_1_inverse_laplace}) into (\ref{u=0}), we can get the expression of $w_0^d(2,t)$. Similarly, using the expressions of $w_0^d(1,t)$ and $w_0^d(2,t)$ we can obtain the expression of $w_0^d(3,t)$.
By applying this iterative algorithm we can identify $w_0^d(n,t)$ for any $n>0$.

Putting the expression of $w_0^d(n,t)$ into the equation (\ref{u>0}) and using the expression of $w_u(k,t,y)$
given in Section \ref{secDickson} allows to calculate the density $w_u^d(n,t)$ for any $u>0$.
We will show later how this algorithm could be used in some examples.

{\bf Remark}
We denote by $ w_u^d(n,t,x)$ the joint density of the number of claims until Parisian ruin time (the corresponding argument
is denoted by $(n)$), the time to Parisian ruin (the corresponding argument
is denoted by $(t)$) and the deficit at Parisian ruin (the corresponding argument
is denoted by $(x)$).
Then similar considerations that gave (\ref{u>0})
gives:
\begin{eqnarray}\label{other joint density}
   &&w_u^d(n,t,x)\\
   &&=\sum_{k=0}^{n-1}\int_{cd}^{\infty}w_{u}(n-k,t-d,y)\frac{(\lambda d)^k}{k!}e^{-\lambda d} f^{k*}(cd-y+x)dy+\sum_{k=1}^{n-1}\int_0^{cd}w_{u}(n-k,t-d,y)\nonumber\\
   &&\times\Big[\frac{(\lambda d)^k}{k!}e^{-\lambda d}f^{k*}(cd-y+x)-\sum_{m=0}^{k-1}\int_{\frac{y}{c}}^dv_y(m,s) \frac{(\lambda (d-s))^{k-m}}{(k-m)!}e^{-\lambda (d-s)}f^{(k-m)*}(c(d-s)+x)ds \Big]dy\nonumber\\
    &&+\sum_{l=1}^{n-1}\sum_{k=0}^{n-l-1}\int_0^{cd}\int_{\frac{y}{c}}^d\int_0^{\max(t-d-s,0)}w_{u}(l,t_1,y)v_y(k,s)w_0^d(n-l-k,t-t_1-s,x)dt_1dsdy.\nonumber\\
 \end{eqnarray}

There is another interesting observation.
Denote the sum of the first term and the second term of formula (\ref{u=0}) by
\begin{eqnarray}
\lefteqn{h(n,t)}\\&:=&\sum_{k=0}^{n-1}\int_{cd}^{\infty}w_{0}(n-k,t-d,y)\frac{(\lambda d)^k}{k!}e^{-\lambda d}dy+\sum_{k=1}^{n-1}\int_0^{cd}w_{0}(n-k,t-d,y)\nonumber\\
   &&\times\Big[\frac{(\lambda d)^k}{k!}e^{-\lambda d}\bar{F}^{k*}(cd-y)-\sum_{m=0}^{k-1}\int_{\frac{y}{c}}^dv_y(m,s) \frac{(\lambda (d-s))^{k-m}}{(k-m)!}e^{-\lambda (d-s)} \bar{F}^{(k-m)*}(c(d-s))ds\Big]dy.\nonumber\\
\end{eqnarray}
Let
\begin{eqnarray*}
\varpi(m,z)&=&\sum_{k=0}^{m-1}\int_0^{cd}w_{0}(m-k,z-s,y)\int_{\frac{y}{c}}^dv_y(k,s)dsdy.
\end{eqnarray*}
Then the equation (\ref{u=0}) can be written more concisely as follows:
\begin{eqnarray}\label{u=0 jian}
    w_0^d(n,t)&=&h(n,t)+\sum_{m=1}^{n-1}\int_0^{t-d}\varpi(m,z)w_0^d(n-m,t-z)dz. \end{eqnarray}
In probability theory, (\ref{u=0 jian}) is known as a (bivariate) renewal equation for the function $w_0^d$.
It is known that the solution of (\ref{u=0 jian}) can be expressed as an infinite series of functions:
\begin{eqnarray}\label{u=0 jiandan}
w_0^d(n,t)&=&\sum_{k=0}^{\infty}\Big(h*\varpi^{k**}\Big)(n,t).
\end{eqnarray}

\subsection{The expression of $p_u^d(n)$}
In this section we will identify $p_u^d(n)$ given in (\ref{pdn}) describing the probability of having $n$ claims up to Parisian ruin time.
Recall that from (\ref{pdn})
\[ p_u(n)=\int_d^\infty w_u^d(n,t)dt.\]
Hence from Theorem \ref {maintheorem} we have the following result.
\begin{theo}
\begin{equation*}
p_u^d(1)=\int_0^{\infty}\lambda e^{-\lambda t}\bar{F}(u+ct+cd)e^{-\lambda d}dt\end{equation*}
and for $n=2,3,\ldots,$
\begin{eqnarray}\label{pudn}
p_u^d(n)&=&\sum_{k=0}^{n-1}\int_{cd}^{\infty}w_{u}(n-k,y)\frac{(\lambda d)^k}{k!}e^{-\lambda d}dy+\sum_{k=1}^{n-1}\int_0^{cd}w_{u}(n-k,y)\nonumber\\
   &&\times\Big[\frac{(\lambda d)^k}{k!}e^{-\lambda d}\bar{F}^{k*}(cd-y)-\sum_{m=0}^{k-1}\int_{\frac{y}{c}}^dv_y(m,s) \frac{(\lambda (d-s))^{k-m}}{(k-m)!}e^{-\lambda (d-s)}\bar{F}^{(k-m)*}(c(d-s))ds\Big]dy\nonumber\\
   &&+\sum_{m=1}^{n-1}\sum_{k=0}^{n-m-1}\int_0^{cd}\int_{\frac{y}{c}}^dw_{u}(m,y)v_y(k,s)p_0^d(n-m-k)dsdy.
\end{eqnarray}
\end{theo}

\section{Examples}

We consider now generic claim size which is exponentially distributed with parameter $\mu$, that is $f(x)=\mu e^{-\mu x}$ for $x>0$.
In this section, we will compute all considered Parisian-type quantities and give some insight on their possible shapes
depending on the choice of parameters.

We start from calculating $p_u^d(n)$.
Let $\lambda=1, \mu=1, c=2$ and $d=2$, we consider different values for the initial surplus, namely: $u=0,1,5,10$.
We show in Table 1 and Figures $3-6$ graphs of $p_0^d(n)$ and $p_u^d(n)$.
\begin{table}[h!]
  \begin{center}
    \begin{tabular}{|c|l|l|l|l|}
      \hline
        \diagbox {n}{$p_u^2(n)$}{u}& 0  & 1 & 5& 10  \\ \hline
       $1$  & 0.0008263 & 0.000303961&$5.56723\times10^{-6} $& $3.75117\times10^{-8}$  \\ \hline
       $2$  & 0.0053243 & 0.00206001 & 0.0000451534& $3.66761\times10^{-7}$  \\ \hline
       $3$  & 0.0129083 & 0.00544101 &0.000156561 &$ 1.62796\times10^{-6} $\\ \hline
       $4$  & 0.0180217 & 0.00848601 &0.00033815&$4.63012\times10^{-6}$   \\ \hline
       $5$  & 0.0179324 & 0.00955855 & 0.000539026&$9.80545\times10^{-6} $ \\ \hline
       $6$  & 0.0146702 & 0.00883697 & 0.000700141& 0.0000168421  \\ \hline
       $7$  & 0.0109439 & 0.00732711 & 0.000790661& 0.0000247947  \\ \hline
       $8$  & 0.0079648 & 0.00577913 &0.00081165&0.0000325081  \\ \hline
       $9$  & 0.0058461 & 0.00448643 &0.000781196& 0.0000390231  \\ \hline
       $10$& 0.0043698 & 0.00348399 &0.000720138& 0.0000437785  \\ \hline
       $11$& 0.0033244 & 0.00272272 &0.000645073& 0.0000466125  \\ \hline
       $12$& 0.0025668 & 0.0013729  &0.000566941 & 0.0000476561  \\ \hline
       $13$& 0.0010261 & 0.00170223 &0.000492033& 0.0000472029  \\ \hline
       $14$& 0.0013668 & 0.00128095 &0.000422024& 0.0000455947  \\ \hline
       $15$& 0.0004854 & 0.00101624 &0.000359118& 0.0000431627 \\ \hline
       $16$& 0.0007728 & 0.00077037 &0.000221917& 0.0000337352  \\ \hline
       $17$& 0.0006557 & 0.000620218 &0.00025548& 0.0000369367 \\ \hline
       $18$& 0.0002086 & 0.000509662 &0.000149864& 0.0000335739  \\ \hline
       $19$& 0.0000977 & 0.000394436 &0.000180647& 0.0000241305 \\ \hline
        \end{tabular}
     \caption*{Tab. 1.   $p_u^2(n)$, $\mu=1$, $\lambda=1, c=2, d=2$.}\label{Ta: first}
   \end{center}
 \end{table}
 \begin{figure}[htbf]
 \centering
 \begin{minipage}[t]{0.45\textwidth}
 \centering
 \includegraphics[scale=0.3]{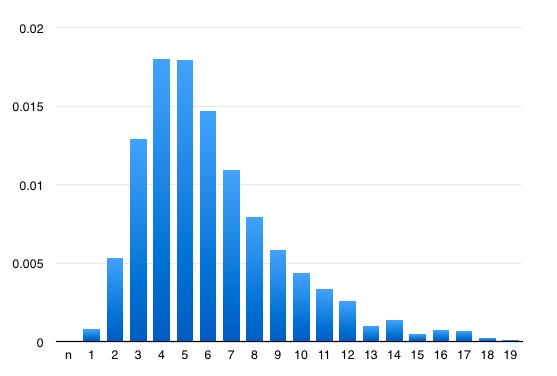}
  \caption*{Fig. 3.   $p_0^2(n)$, $\mu=1$, $\lambda=1, c=2, d=2$.}
  \end{minipage}
  \begin{minipage}[t]{0.45\textwidth}
   \centering
  \includegraphics[scale=0.3]{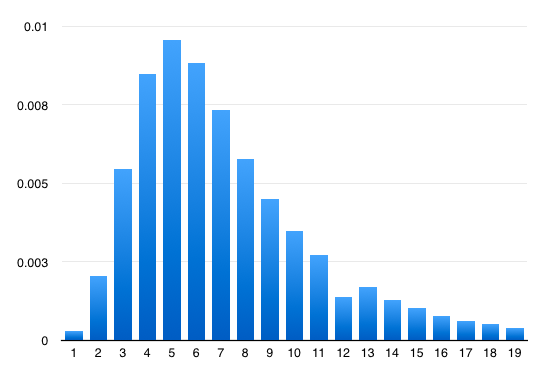}
  \caption*{Fig. 4.  $p_1^2(n)$, $\mu=1$, $\lambda=1, c=2, d=2$.}
  \end{minipage}
 \end{figure}

 \begin{figure}[htbf]
 \centering
 \begin{minipage}[t]{0.45\textwidth}
 \centering
 \includegraphics[scale=0.3]{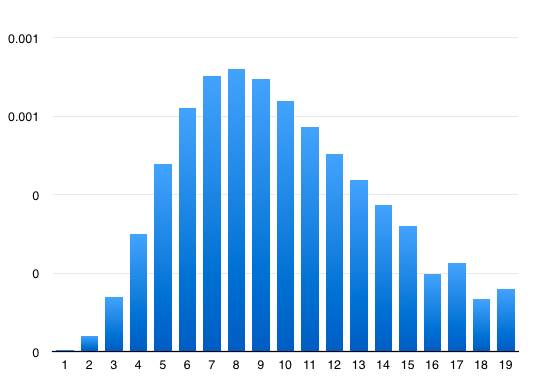}
  \caption*{Fig. 5.  $p_5^2(n)$, $\mu=1$, $\lambda=1, c=2, d=2$.}
  \end{minipage}
  \begin{minipage}[t]{0.45\textwidth}
   \centering
  \includegraphics[scale=0.3]{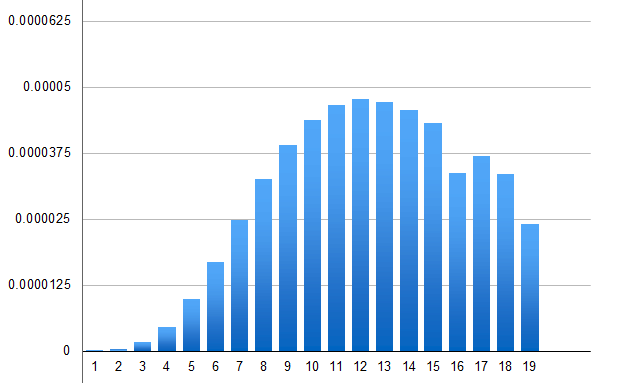}
  \caption*{Fig. 6.   $p_{10}^2(n)$, $\mu=1$, $\lambda=1, c=2, d=2$.}  \end{minipage}
 \end{figure}

We noticed the following observations. The distributions are quite cumulated around the mode of $p_u^d(n)$ as a function $n$.
Moreover, for fixed $u$, these probabilities first increase then decrease when the number of claims gets bigger.
It seems the tails of this probability functions is surprisingly thick. In fact, it seems that larger $u$ produces thicker tail.

Now we will focus on more complex density $w_u^d(n,t)$.
We will find it using Theorem \ref{maintheorem}.

From (\ref{w u d 1}) it follows that:
\begin{eqnarray}\label{exp w u d 1}
w_u^d(1,t)&=&\lambda e^{-(\lambda+\mu c) t-\mu u}.
 \end{eqnarray}
We will calculate now $w_0^d(n,t)$ for $n>1$ and $u=0$.
Using \ref{expv} and (\ref{w0kt}) from we can derive the expression for $w_0^d(2,t)$:
\begin{eqnarray}\label{w0d2t}
w_0^d(2,t)=\left\{
  \begin{array}{lll}
&\lambda^2e^{-(\lambda+\mu c) t}(\frac{1}{2}\mu ct^2+d-\frac{1}{2}\mu cd^2),&t>2d,\\
&\lambda^2e^{-(\lambda+\mu c) t}(\mu c(t-d)^2+d+\frac{1}{2}\mu cd^2),&d<t\leq2d.
  \end{array}
\right.
\end{eqnarray}
Similarly, using the expressions of $w_0^d(1,t)$ and $w_0^d(2,t)$ we can obtain the expression of $w_0^d(3,t)$:
\begin{eqnarray}\label{w0d2t}
w_0^d(3,t)=\left\{
  \begin{array}{lll}
&\lambda^3e^{-(\lambda+\mu c) t}f_{31}(t),&t>3d,\\
&\lambda^3e^{-(\lambda+\mu c) t}f_{32}(t),&2d<t\leq3d\\
&\lambda^3e^{-(\lambda+\mu c) t}f_{33}(t),&d<t\leq2d,
  \end{array}
\right.
\end{eqnarray}
where
\begin{eqnarray*}
&&f_{31}(t)=\frac{1}{24}[2c^2\mu^2t^4+(12cd\mu-6c^2d^2\mu^2)t^2+4c^2d^3\mu^2t +12d^2-16cd^3\mu+5c^2d^4\mu^2],\\
&&f_{32}(t)=\frac{1}{24}[c^2\mu^2t^4+12c^2d\mu^2t^3+(12cd\mu-60c^2d^2\mu^2)t^2+112c^2d^3\mu^2t+12d^2-16cd^3\mu-76c^2d^4\mu^2],\\
&&f_{33}(t)=\frac{1}{12}[3c^2\mu^2t^4-12c^2d\mu^2t^3+(12cd\mu+24c^2d^2\mu^2)t^2-(24cd^2\mu+24c^2 d^3\mu^2)t+6d^2+16cd^3\mu\\
&&\ \ \ \ \ \ \ \ \ \ \ \ +10c^2d^4\mu^2].
  \end{eqnarray*}
This iterative algorithm can produce $w_0^d(n,t)$ for any $n>0$ and then by Theorem \ref{maintheorem} we can identify all $w_u^d(n,t)$.
Unfortunately, the computation process takes long time and the expression for $w_u^d(n,t)$  gets very complicated quite quickly.
We suggest another numerical algorithm instead.
We change the integration in the third increment of (\ref{u=0}) and (\ref{u>0}) into the
summation using the rectangular method of approximating definite integrals.
Of course taking the step of the summation tending to $0$ will give right expression.

Let $\lambda=1, \mu=1, c=1.2, d=2$. Here we divide up the interval $[0,10]$ into $100$ equal subintervals.
Each has length  $\Delta t = \dfrac{1}{10}$ . We will evaluate the function $w_0^d(n,t)$ at the right-hand endpoints of these subinterval.
Figure 7 shows that the approximation is very good. We noticed that the time of calculating $w_0^d(n,t)$ is now much shorter.
(Solid line denote the exactly values of the density function $w_0^d(n,t)$ and dotted line denote the approximate values.)

 \begin{figure}[!hbp]
 \centering
 \begin{minipage}[t]{0.45\textwidth}
 \centering
 \includegraphics[scale=0.5]{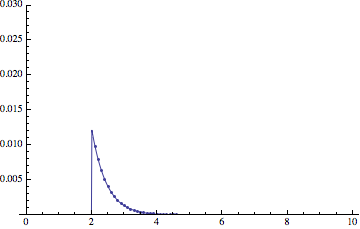}
 \caption*{ $w_0^d(1,t)$, $\mu=1$, $\lambda=1, c=1.2, d=2$.}     \end{minipage}
  \begin{minipage}[t]{0.45\textwidth}
   \centering
  \includegraphics[scale=0.5]{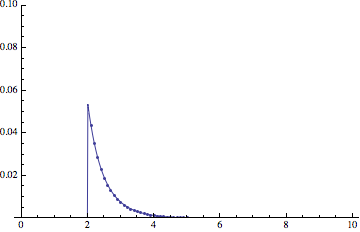}
  \caption*{$w_0^d(2,t)$, $\mu=1$, $\lambda=1, c=1.2, d=2$.}     \end{minipage}
 \end{figure}

  \begin{figure}[!hbp]
 \centering
 \begin{minipage}[t]{0.45\textwidth}
 \centering
 \includegraphics[scale=0.5]{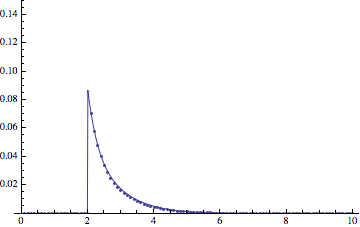}
 \caption*{$w_0^d(3,t)$, $\mu=1$, $\lambda=1, c=1.2, d=2$.}    \end{minipage}
  \begin{minipage}[t]{0.45\textwidth}
   \centering
 \end{minipage}
 \caption*{Fig. 7.  The exact and approximate values of $w_0^d(n,t)$. }
  \end{figure}

Figure 8 shows the graphs for $w_0^d(n,t)$ for $n=1,2,3,4,5,6$ and Figure 9 shows the graphs for $w_u^d(n,t)$ when $d=2$, $u=2$.
Tables $2-3$ give the values of $p_0^d(n,t)$ and $p_u^d(n,t)$ for $n=1,2,3,4,5,6$, $t=1,2,3,4,5,6,7$, respectively.

  \begin{figure}[!hbp]
 \centering
 \begin{minipage}[t]{0.45\textwidth}
 \centering
 \includegraphics[scale=0.5]{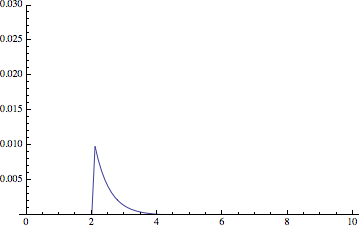}
 \caption*{$w_0^d(1,t)$, $\mu=1$, $\lambda=1, c=1.2, d=2$.}   \end{minipage}
  \begin{minipage}[t]{0.45\textwidth}
   \centering
  \includegraphics[scale=0.5]{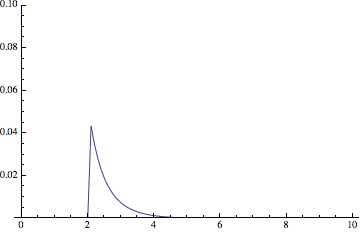}
  \caption*{$w_0^d(2,t)$, $\mu=1$, $\lambda=1, c=1.2, d=2$.}  \end{minipage}
 \end{figure}

 \begin{figure}[!hbp]
 \centering
 \begin{minipage}[t]{0.45\textwidth}
 \centering
 \includegraphics[scale=0.5]{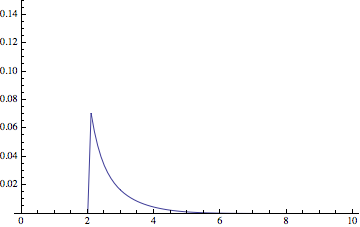}
 \caption*{$w_0^d(3,t)$, $\mu=1$, $\lambda=1, c=1.2, d=2$.}   \end{minipage}
  \begin{minipage}[t]{0.45\textwidth}
   \centering
  \includegraphics[scale=0.5]{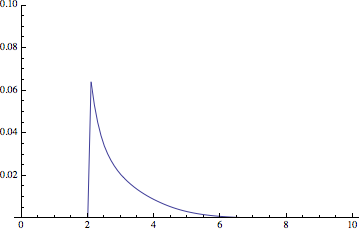}
  \caption*{$w_0^d(4,t)$,$\mu=1$, $\lambda=1, c=1.2, d=2$.}  \end{minipage}
 \end{figure}

\begin{figure}[!hbp]
 \centering
 \begin{minipage}[t]{0.45\textwidth}
 \centering
 \includegraphics[scale=0.5]{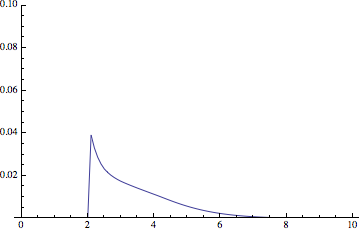}
 \caption*{$w_0^d(5,t)$, $\mu=1$, $\lambda=1, c=1.2, d=2$.}   \end{minipage}
  \begin{minipage}[t]{0.45\textwidth}
  \centering
  \includegraphics[scale=0.5]{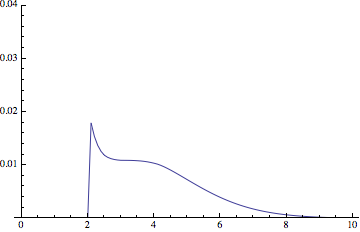}
 \caption*{$w_0^d(6,t)$, $\mu=1$, $\lambda=1, c=1.2, d=2$.}  \end{minipage}
  \caption*{Fig. 8.  Graphs for $w_0^d(n,t)$, $\mu=1$, $\lambda=1, c=1.2, d=2$.}
  \end{figure}

\begin{figure}[!hbp]
\centering
 \begin{minipage}[t]{0.45\textwidth}
  \centering
  \includegraphics[scale=0.5]{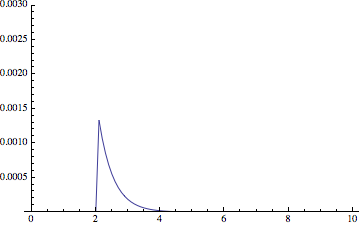}
 \caption*{$w_u^d(1,t)$, $\mu=1$, $\lambda=1, c=1.2, d=2, u=2$.}  \end{minipage}
\begin{minipage}[t]{0.45\textwidth}
 \centering
 \includegraphics[scale=0.5]{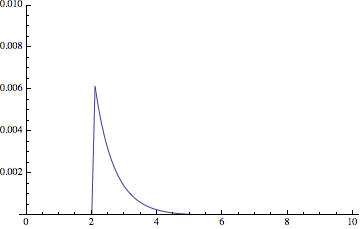}
 \caption*{$w_u^d(2,t)$, $\mu=1$, $\lambda=1, c=1.2, d=2, u=2$.}   \end{minipage}
 \end{figure}

 \begin{figure}[!hbp]
 \centering
  \begin{minipage}[t]{0.45\textwidth}
   \centering
  \includegraphics[scale=0.5]{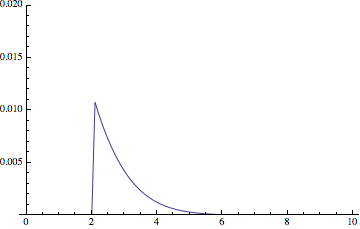}
  \caption*{$w_u^2(3,t)$, $\mu=1$, $\lambda=1, c=1.2, d=2, u=2$.}  \end{minipage}
   \begin{minipage}[t]{0.45\textwidth}
 \centering
 \includegraphics[scale=0.5]{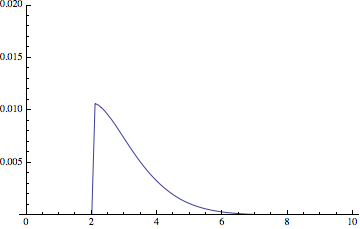}
 \caption*{$w_u^2(4,t)$, $\mu=1$, $\lambda=1, c=1.2, d=2, u=2$.}   \end{minipage}
  \end{figure}

\begin{figure}[!hbp]
 \centering
 \begin{minipage}[t]{0.45\textwidth}
   \centering
  \includegraphics[scale=0.5]{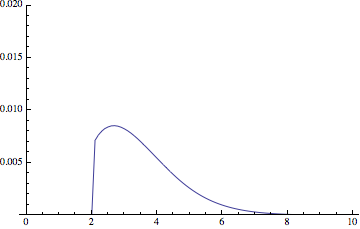}
  \caption*{ $w_u^2(5,t)$, $\mu=1$, $\lambda=1, c=1.2, d=2, u=2$.}  \end{minipage}
   \begin{minipage}[t]{0.45\textwidth}
 \centering
 \includegraphics[scale=0.5]{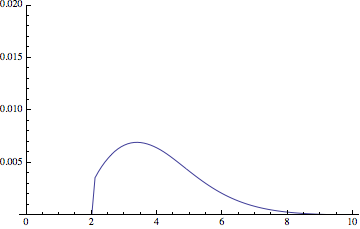}
 \caption*{ $w_u^d(6,t)$,$\mu=1$, $\lambda=1, c=1.2, d=2, u=2$.}   \end{minipage}
 \caption*{Fig. 9.  Graphs for $w_u^d(n,t)$, $\mu=1$, $\lambda=1, c=1.2, d=2, u=2$.}
  \end{figure}

 \begin{table}[!hbp]
  \begin{center}
    \begin{tabular}{|c|c|c|c|c|c|c|c|}
      \hline
       \diagbox{n} {$p_0^d(n,t)$} {t}&$t=1$&$t=2$&$t=3$&$t=4$&$t=5$&$t=6$&$t=7$ \\ \hline
       $n=1$  &0&0&0.0044364&0.00492798&0.00498245&0.00498848&0.00498915\\ \hline
       $n=2$  &0&0&0.0204411&0.0236396&0.0242160&0.0243135&0.0243288\\ \hline
       $n=3$  &0&0&0.0360174&0.0445814&0.0469080&0.0474601&0.0475778\\ \hline
       $n=4$  &0&0&0.0360946&0.0494625&0.0548542&0.0566279&0.0571321\\ \hline
       $n=5$  &0&0&0.0246765&0.0386667&0.0468287&0.0505039&0.0518839\\ \hline
       $n=6$  &0&0&0.0127279&0.0234754&0.0323266&0.0377145&0.0403638\\ \hline
       $n=7$  &0&0&0.00527283&0.0117139& 0.019034&0.0249818&0.0287784\\ \hline
       $n=8$  &0&0&0.00182873&0.00497309&0.00980279&0.0149746&0.0192228\\ \hline
  \end{tabular}
     \caption*{Tab. 2.   $p_0^d(n,t)$, $\mu=1$, $\lambda=1, c=1.2, d=2$.}\label{Ta: first}
   \end{center}
 \end{table}

 \begin{table}[!hbp]
  \begin{center}
    \begin{tabular}{|c|c|c|c|c|c|c|c|}
      \hline
 \diagbox{n} {$p_u^d(n,t)$} {t}& $t=1$  & $t=2$ & $t=3$ & $t=4$ &$t=5$&$t=6$&$t=7$\\ \hline
       $n=1$  &0&0&0.000600403&0.000666929&0.000674301&0.000675117&0.000675208\\ \hline
       $n=2$  &0&0&0.00322483&0.00384155&0.00395467&0.00397339&0.00397626\\ \hline
       $n=3$  &0&0&0.00711875&0.00943587&0.0100615&0.0102047&0.0102339\\ \hline
       $n=4$  &0&0&0.00920644&0.0141209&0.0160422&0.0166421&0.016803 \\ \hline
       $n=5$  &0&0&0.00816603&0.0149827&0.0187639&0.0203559&0.0209108\\ \hline
       $n=6$  &0&0&0.00542406&0.0122129&0.0174604&0.0204087&0.0217345\\ \hline
          \end{tabular}
     \caption*{Tab. 3.   $p_u^d(n,t)$, $\mu=1$, $\lambda=1, c=1.2, d=2, u=2$.}\label{Ta: first}
   \end{center}
 \end{table}

\noindent{\bf{Acknowledgements}}\\
This research is support by the National Natural Science Foundation of China (Grant No. 11271164 and 11371020) and by the FP7 Grant PIRSES-GA-2012-318984.
Zbigniew Palmowski is supported by the Ministry of Science and Higher Education of Poland under the grant
2013/09/B/HS4/01496.

\bibliographystyle{plain}
\bibliography{refbib}
\end{document}